\documentclass[12pt,draftcls,onecolumn]{IEEEtran}  
\usepackage[dvips]{epsfig}
\usepackage{color}
\usepackage{graphicx}
\usepackage[tight,footnotesize]{subfigure}
\usepackage{authblk, amsmath, amssymb}

\newtheorem{theorem}            {Theorem}[section]

\begin{document}
\title{Optimal stopping for the predictive maintenance of a structure subject to  corrosion}

\author[1,2]{Beno\^{\i}te de Saporta} 
\author[1]{Fran\c cois Dufour}
\author[1]{Huilong Zhang}
\author[3]{Charles Elegbede}

\affil[1]{Universit\'e de Bordeaux, IMB, CNRS UMR~5251 INRIA~Bordeaux~Sud~Ouest team CQFD}
\affil[2]{Universit\'e de Bordeaux, GREThA, CNRS UMR~5113}
\affil[3]{Astrium}

\maketitle
\begin{abstract}
We present a numerical method to compute the optimal maintenance time for a 
complex dynamic system applied to an example of maintenance of a metallic structure subject to corrosion. 
An arbitrarily early intervention may be uselessly costly, but a 
late one may lead to a partial/complete failure of the system, which has to 
be avoided. One must therefore find a balance between these too simple 
maintenance policies. To achieve this aim, we model the system by a stochastic 
hybrid process. The maintenance problem thus corresponds to an optimal stopping problem. 
We propose a numerical method to solve the optimal stopping problem and optimize the maintenance 
time for this kind of processes.
\end{abstract}


\section{Introduction}
A complex system is inherently sensitive to failures of its components. 
We must therefore determine maintenance policies in order to maintain an 
acceptable operating condition. The optimization of maintenance is a very 
important problem in the analysis of complex systems. It determines when 
maintenance tasks should be performed on the system. These intervention dates should be chosen to 
optimize a cost function, that is to say, maximize a performance function or, 
similarly, to minimize a loss function. Moreover, this optimization must take
 into account the random nature of failures and random evolution 
and dynamics of the system. Theoretical study of the optimization of 
maintenance is also a crucial step in the process of optimization of 
conception and study of the life service of the system before the first maintenance. 

\bigskip

We consider here an example of maintenance related to an aluminum metallic 
structure subject to corrosion. This example was provided by Astrium. It concerns a small structure within a strategic ballistic missile. The missile is stored successively in a workshop, in a nuclear submarine missile 
launcher in operation or in the submarine in dry-dock. These various environments are more or less corrosive and the structure is inspected with a given periodicity. It is 
made to have potentially large storage durations. The requirement 
for security is very strong. The mechanical stress exerted 
on the structure depends in part on its thickness. A loss of thickness will cause an 
over-constraint and therefore increase a risk of rupture. It is thus crucial to 
control the evolution of the thickness of the structure over time, and to 
intervene before the failure. 

\bigskip

The only maintenance operation we consider here is the complete replacement of the structure. We do not allow partial repairs. Mathematically, this problem of preventive maintenance 
corresponds to a stochastic optimal stopping problem as explained by example 
in the book of Aven and Jensen \cite{aven99a}. It is a difficult problem, because on the one hand, the structure spends random times in each environment, and on the other hand, the corrosiveness of each  environment is also supposed to be random within a given range. In addition, we search for an optimal maintenance date adapted to the particular history of each structure, and not an average one. We also want to be able to update the predicted maintenance date given the past history of the corrosion process.

\bigskip

To solve this maintenance problem, we propose to model this system by a 
piecewise-de\-ter\-mi\-nis\-tic Markov process 
 (PDMP). PDMP's are a class of stochastic hybrid processes that have been 
introduced by Davis \cite{Davis93a} in the 80's. These processes have two 
components: a Euclidean component that represents the physical system 
(e.g. temperature, pressure, thickness loss) and a discrete component that 
describes its regime of operation and/or its environment. Starting from a state $x$ and 
mode $m$ at the initial time, the process follows a deterministic trajectory given 
by the laws of physics until a jump time that can be either 
random (e.g. it corresponds to a component failure or a change of environment) 
or deterministic (when a magnitude reaches a certain physical threshold, 
for example the pressure reaches a critical value that triggers a valve). 
The process restarts from a new state and a new mode of operation, and so on. 
This defines a Markov process. Such processes can naturally take into account 
the dynamic and uncertain aspects of the evolution of the system. 
A subclass of these processes has been introduced by Devooght 
\cite{devooght97a} for an application in the nuclear field. The general 
model has been introduced in dynamic reliability by Dutuit and Dufour 
\cite{dufour02a}. 

\bigskip

The  theoretical problem of optimal stopping for PDMP's is well understood, see e.g.
Gugerli \cite{gugerli86a}. However, there are surprisingly few works in the literature presenting practical algorithms to compute the optimal cost and optimal stopping time. To our best knowledge only 
Costa and Davis 
\cite{costa88a} have presented an algorithm for calculating these quantities 
for PDMP's. Yet, as illustrated above, it is 
crucial to have an efficient numerical tool to compute the  optimal maintenance 
time in practical cases. The purpose of this paper is to adapt the general algorithm recently proposed by the authors in \cite{saporta10a} to this special case of maintenance and show its high practical power. More precisely, we present a method to compute the optimal cost as well as a quasi optimal stopping rule, that is the date when the maintenance should be performed. As a byproduct of our procedure, we also obtain the distribution of the optimal maintenance dates and can compute dates such that the probability to perform a maintenance before this date is below a prescribed threshold.

\bigskip

The remainder of this paper is organized as follows. In section~\ref{modeling}, 
we present the example of corrosion of the metallic structure that  we are 
interested in with more details as well as the framework of PDMP's. 
In section~\ref{opt stop}, we briefly recall the formulation of the optimal 
stopping problem for PDMP's and its theoretical solution. In section~\ref{num}, 
we detail the four main steps of algorithm. In section~\ref{results} we present the numerical 
results obtained on the example of corrosion. Finally, in section~\ref{ccl}, 
we present a conclusion and perspectives.

\section{Modeling}
\label{modeling}

Throughout this paper, our approach will be illustrated on an example of 
maintenance of a metallic structure subject to corrosion. This example was proposed by Astrium.
As explained in the introduction, it is a small homogeneous aluminum structure within a strategic ballistic missile.
The missile is stored for potentially long times in more or less corrosive environments.
The mechanical stress exerted on the structure depends in part on its thickness. A loss of thickness will cause an 
over-constraint and therefore increase a risk of rupture. It is thus crucial to 
control the evolution of the thickness of the structure over time, and to 
intervene before the failure. 

\bigskip

Let us describe more precisely the usage profile of the missile. Its is stored successively in three different environments, the workshop, the submarine in operation and the submarine in dry-dock. This is 
because the structure must be equipped and used in a given order. Then it goes back to the workshop and so on. The missile stays in each environment during a random duration with exponential distribution. Its parameter depends on the environment. At the beginning of its service time, the structure is treated against corrosion. The period of effectiveness of this protection is also random, with a Weibull distribution. The thickness loss only begins when this initial protection is gone. The degradation law for the thickness loss then depends on the environment through two parameters, a deterministic transition period and a random corrosion rate uniformly distributed within a given range. Typically, the workshop and dry-dock are the more corrosive environments. The randomness of the corrosion rate accounts for small variations and uncertainties in the corrosiveness of each environment. 

\bigskip

We model this degradation process by a $3$-dimensional PDMP ($X_t$)  with 3 modes corresponding to the three different environment. Before giving the detailed parameters of this process, we shortly present general PDMP's.
 
\subsection{Definition of piecewise-deterministic Markov processes}

Piecewise-deterministic Markov processes (PDMP's) are a general class of hybrid processes. Let $M$ be the finite set of the possible modes of the system. In our example, the modes correspond to the various environments. For all mode $m$ in $M$, 
let $E_m$ an open subset in $\mathbb{R}^d$.  A PDMP is 
defined from three local characteristics $(\Phi, \lambda, Q)$ where
\begin{itemize} 
\item the flow $\Phi : M\times \mathbb{R}^d\times \mathbb{R} \rightarrow \mathbb{R}^d$ is continuous 
and for all $ s, t \geq 0$,  one has $\Phi(\cdot, \cdot,t+s) = \Phi(\Phi(\cdot, \cdot,s),t)$. 
It describes 
the deterministic trajectory of the process between jumps. For all $(m, x)$ 
in $ M \times E_m$, we set
\begin{displaymath} 
\displaystyle 
t^\ast(m,x) = \inf \{t>0 : \Phi(m,x,t) \in \partial E_m\},
\end{displaymath} 
the time to reach the  boundary of the domain starting from $x$ in mode $m$.
\item the jump intensity $\lambda$ characterizes the frequency of jumps. For 
all $(m, x)$ in $M \times E_m$, and $t \leq t^\ast (m, x)$, we set
\begin{displaymath} 
\displaystyle 
\Lambda(m,x,t) = \int_0^t \lambda (\Phi(m,x,s))\,ds.
\end{displaymath}  
\item the Markov kernel $Q$ represents the transition measure of the process 
and allows to select the new location after each jump.
\end{itemize} 

\begin{figure}[h]
\centering
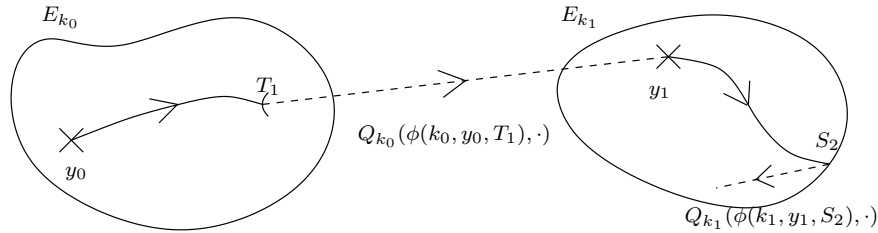
\caption{An exemple of path for a PDMP until the second jump. The first jump is 
random. The second jump is deterministic because the process has reached the 
boundary of the domain.}
\label{figure_1}
\end{figure}

The trajectory $X_t = (m_t, x_t)$ of the process can then be defined 
iteratively. We start with an initial point $X_0 = (k_0, y_0)$ with $ k_0\in M$
 and $y_0\in E_{k_0}$. The first jump time $T_1$ is determined by
\begin{displaymath} 
\displaystyle 
\mathbb{P}_{(k_0,y_0)}(T_1>t) = \left \{ 
\begin{array} {lcl}
e ^{-\Lambda(k_0,y_0,t)} &\textrm{if} &~t<t^\ast(k_0,y_0),\\
0 &\textrm{if}  &~t\geq t^\ast(k_0,y_0).
\end{array} 
\right .
\end{displaymath} 
On the interval $[0, T_1)$, the process follows the  deterministic trajectory 
$m_t = k_0$ and $x_t = \Phi(k_0, y_0, t)$. At the random time $T_1$, 
a jump occurs. Note that a jump can be either a discontinuity in the Euclidean variable $x_t$ or a change of mode. The process restarts at a new mode and/or position $X_{T_1}=(k_1,y_1)$, according 
to distribution $Q_{k_0}(\Phi(k_0, y_0, T_1), \cdot)$. We then select in a similar way an
inter jump time  $T_2- T_1$, and in the interval $[T_1, T_2)$ the process 
follows the path $m_t=k_1$ and $x_t = \Phi(k_1, y_1, t - T_1)$. 
Thereby, iteratively, a PDMP is constructed, see Figure \ref{figure_1} for an 
illustration. 
Let $Z_0 = X_0$, and for $n \geq 1, ~Z_n = X_{T_N}$, location and mode of the 
process after each jump. Let $S_0 = 0$, $S_1=T_1$ and for $n \geq 2$, 
$S_n = T_n-T_{n-1}$ the inter-jump times between two consecutive jumps, 
then $(Z_n, S_n)$  is a Markov chain, which is the only source of 
randomness of the PDMP and contains all information on its random part.
 Indeed, if one knows the jump times and the positions after each jump, we can 
reconstruct the deterministic part of the trajectory between jumps. It is a 
very important property of PDMP's that is at the basis of our numerical procedure.

\subsection{Example of corrosion of metallic structure}
\label{section_example_corrosion}
We can  now turn back to our example of corrosion of structure and give the characteristics of the PDMP modeling the thickness loss. The finite set of modes is $M=\{1,2,3\}$, where mode $1$ corresponds to the workshop environment, mode $2$ to the submarine in operation and mode $3$ to the dry-dock. Although the thickness loss is a one-dimensional process, one needs a three dimensional PDMP to model its evolution, because it must also take into account all the sources of randomness, that is the duration of the initial protection and the corrosion rate in each  environment. The corrosion process ($X_t$) is defined by:
\begin{displaymath} 
\displaystyle 
X_t = (m_t,d_t,\gamma_t,\rho_t)\in\{1,2,3\}\times\mathbb{R}_+\times\mathbb{R}_+\times\mathbb{R}_+,
\end{displaymath} 
where $m_t$ is the environment at time $t$, $d_t$ is the thickness loss at time $t$, $\gamma_t$ is the remainder of the initial protection at time $t$ and $\rho_t$ is the corrosion rate of the current environment at time $t$.

\bigskip

Originally, at time $0$, one has $X_0=(1,0,\gamma_0,\rho_0)$, which means that the missile is in the workshop and the structure has not started corroding yet. The original protection $\gamma_0$ is drawn according to a Weibull distribution function
\begin{displaymath} 
\displaystyle 
F(t) = 1 - \exp \left (-\left (\frac{t}{\beta}\right )^\alpha\right )
\end{displaymath} 
with $\alpha=2.5$ and $\beta=11800$ hours$^{-1}$. The corrosion rate in the workshop is drawn according to a uniform distribution on $[10^{-6}, 10^{-5}]$ mm/hour. The time $T_1$ spent in the workshop is drawn according to an exponential distribution with parameter $\lambda_1 = 17520$ hour$^{-1}$. At time $t$ between time $0$ and time $T_1$, the remainder of the protection is simply $\gamma_t=\max\{0,\gamma_0-t\}$, $\rho_t$ is constant equal to $\rho_0$ and the thickness loss $d_t$ is given by
\begin{equation} 
d_t = \left \{
\begin{array} {lcl}
0 &\textrm{if}& ~ t\leq \gamma_0, \\
\displaystyle \rho_0 \left (
t-(\gamma_0+\eta_1)+\eta_1\exp\left (-\frac{t-\gamma_0}{\eta_1}\right)
\right ) &\textrm{if}& ~ t>\gamma_0,
\end{array} 
\right .
\label{equa_corrosion}
\end{equation} 
where $\eta_1=30000$ hours.

\bigskip

At time $T_1$, a \emph{jump} occurs, which means there is a change of environment and a new corrosion rate is drawn for the new environment. The other two components of the process $(X_t)$ modeling the remainder of the protection $\gamma_t$ and the thickness loss $d_t$ naturally evolve continuously. Therefore, one has $m_{T_1} = 2$, $\gamma_{T_1} = 0$ if $\gamma_0<T_1$, $\gamma_{T_1} = \gamma_0-T_1$ 
otherwise~; that is to say that once the initial protection is gone, 
it has no effect any longer, $\rho_{T_1}$ is drawn according to a uniform distribution on $[10^{-7}, 10^{-6}]$ mm/hour.
The process continues to evolve in the same way until the next change of environment occurring at time $T_2$. Between $T_1$ and $T_2$, just replace $\rho_0$ by $\rho_{T_1}$, $\gamma_0$ by $\gamma_{T_1}$, $\eta_1$ by $\eta_2=200000$ hours and $t$ by $t-T_1$ in equation~(\ref{equa_corrosion}). The process visits successively the 
3 environments always in the same order 1, 2 and 3 and then returns to the environment 1. 
. The time spent in the environment $i$ is a random variable exponentially distributed with 
parameters $\lambda_i$ with $\lambda_1= 17520$ hours$^{-1}$, $\lambda_2 = 131400$ hours$^{-1}$
and $\lambda_3 = 8760$ hours$^{-1}$. The thickness loss evolves continuously 
according to equation (\ref{equa_corrosion}) with suitably changed parameters. The period of transition in the 
mode $i$ is 
$\eta_i$ with $\eta_1 = 30000$~hours, $\eta_2 = 200000$~hours and  $\eta_3 = 40000$~hours. 
The corrosion rate $\rho_i$ expressed in mm per hour is drawn at each change of 
environments. In environments 1 and 3, it follows a uniform distribution
on $[10^{-6}, 10^{-5}]$ and in environment 2, it follows a uniform distribution
on $[10^{-7}, 10^{-6}]$. 
\begin{figure}[ht]
\centering
\subfigure[One trajectory]
{\includegraphics[width=0.48\linewidth]{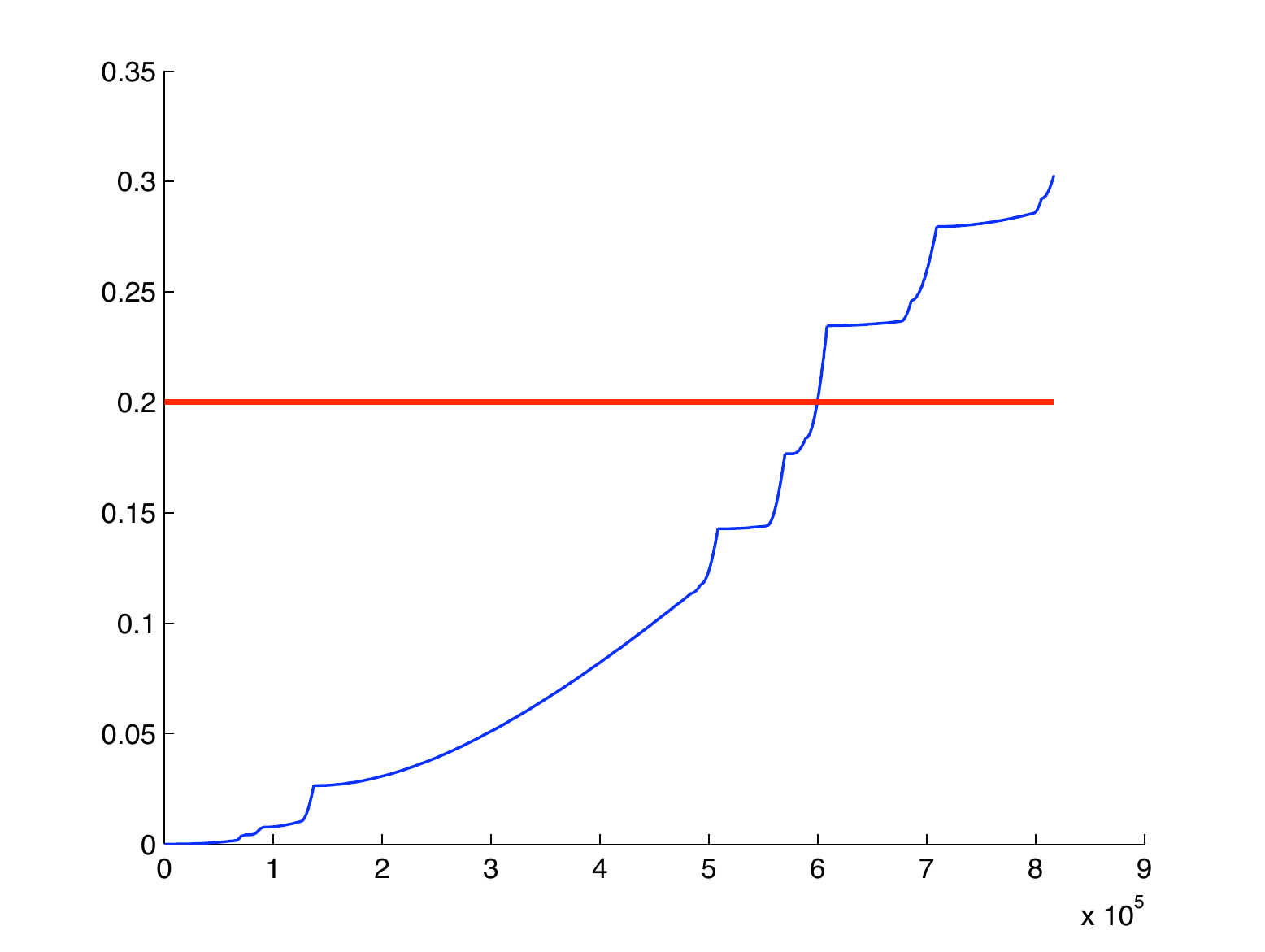}} 
\subfigure[100 trajectories]
{\includegraphics[width=0.48\linewidth]{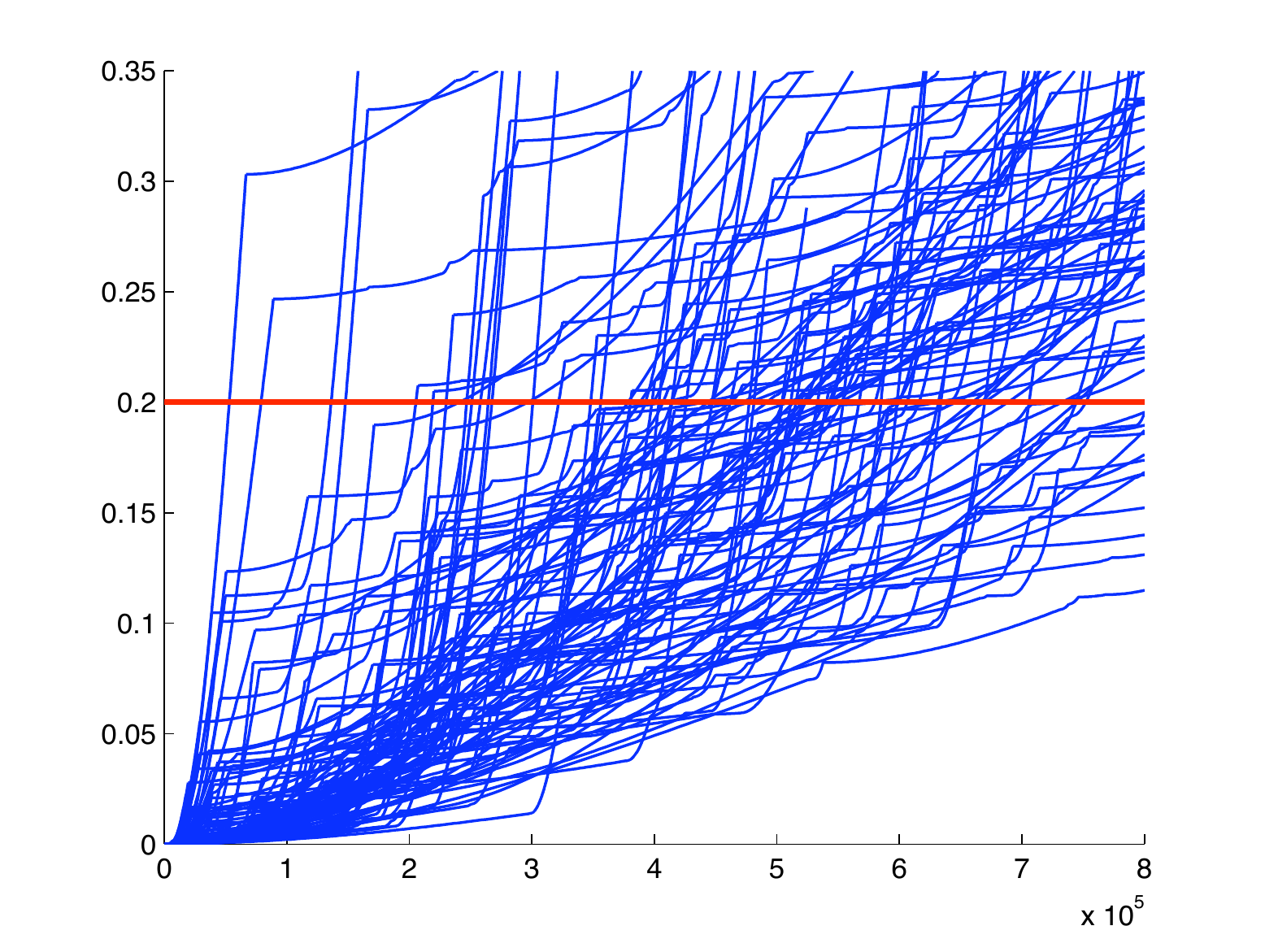}}
\caption{Examples of trajectories of thickness loss over time.}
\label{figure_2}
\end{figure}
Figure \ref{figure_2} shows examples of simulated 
trajectories of the thickness loss. The slope changes correspond 
to changes of environment. The observed dispersion is characteristic of the 
random nature of the phenomenon. Note that the various physical parameters were given by Astrium and will not be discussed here.

\bigskip

The missile is inspected and the thickness loss of the structure under study is measured at each change of environment. Note that the structure is small enough for only one measurement point to be significant. The structure is considered unusable if the loss of thickness reaches $0.2$mm. 
The optimal maintenance time must therefore occur before reaching this 
critical threshold, which could cause the collapse of the 
structure, but not too soon which would be unnecessarily expensive. It should also only use the available measurements of the thickness loss.

\section{Optimal stopping problem}
\label{opt stop}

We now briefly formulate the general mathematical problem of optimal stopping corresponding to our maintenance problem. Let 
$z = (k_0, y_0)$ be the starting point of the PDMP $(X_t)$. Let $\mathcal{M}_N$ be the 
set of all stopping times $T$ for the natural filtration of the PDMP ($X_t$) satisfying 
$T \leq T_N$ that is to say that the intervention takes place before the $N$th jump of process. 
The $N$th jump represents the horizon of our maintenance problem, that is to 
say that we impose to intervene no later than $N$th change of environment. The choice of $N$ is discussed below. Let $g$
 be the cost function to optimize. Here, $g$ is a reward function that we
want to maximize. The optimization problem to solve is the following
\begin{displaymath} 
\displaystyle 
v(z) = \sup_{\tau \in M_N}E_z\left [g(X_\tau)\right ].
\end{displaymath} 
The function $v$ is called the {\em value function} of the problem and represents 
the maximum performance that can be achieved. Solving the optimal stopping 
problem is firstly to calculate the value function, and secondly to find a 
stopping time $\tau$ that achieves this maximum. This stopping time is 
important from the application point of view since it corresponds to the 
optimum time for maintenance. 
In general, such an optimal stopping time does not exist. We then define 
$\epsilon$-optimal stopping times as achieving optimal value minus $\epsilon$, 
i.e. $v(z)-\epsilon$.

\bigskip

Under fairly weak regularity conditions, Gugerli has shown in \cite{gugerli86a} 
that the value function $v$ can be calculated iteratively as follows. 
Let $v_N=g$ be the reward function, and we iterate an operator $L$ backwards. 
The function $v_0$ thus obtained is equal to the value function $v$.
\begin{displaymath} 
\left \{
\begin{array} {lcl}
v_N & = & g, \\
v_k & = & L(v_{k+1},g), \quad0\leq k\leq N-1.
\end{array} 
\right .
\end{displaymath} 
The operator $L$ is a complex operator which involves a 
continuous maximization, conditional expectations and indicator functions, 
even if the cost function $g$ is very regular.
\begin{equation} \label{def L}
\begin{array} {rl}
\displaystyle L(w,g)(z) 
\equiv & \displaystyle \sup_{u\leq t^\ast(z) }
\left \{\displaystyle 
E \left [w(Z_1)1_{S_1<u\wedge t^\ast(z)} 
+ g(\Phi(z,u))1_{S_1\geq u \wedge t^\ast(z)}|Z_0=z \right ] 
\right \}\\
 & \vee E \left [ w(Z_1)|Z_0=z \right ].
\end{array} 
\end{equation} 
However, we can see that this operator depends only on the discrete time 
Markov chain $(Z_n, S_n)$. Gugerli also proposes an iterative construction of $\epsilon$-optimal 
stopping times, which is a bit too tedious and technical to be described here, 
see \cite{gugerli86a} for details.

\bigskip

\begin{figure}[ht]
\centering
\subfigure
{\includegraphics[width=0.50\linewidth]{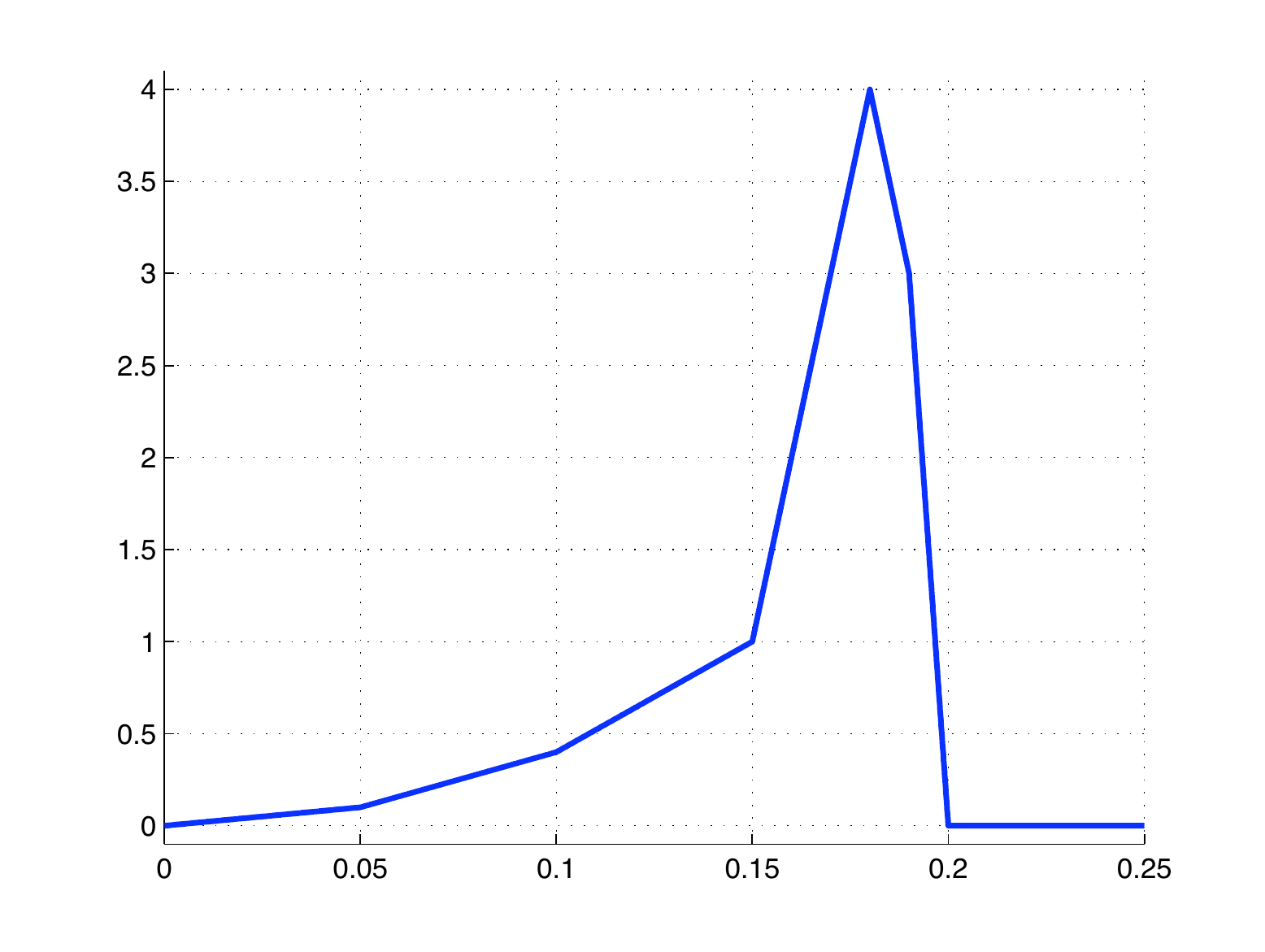}} 
\subfigure
{\includegraphics[width=0.48\linewidth]{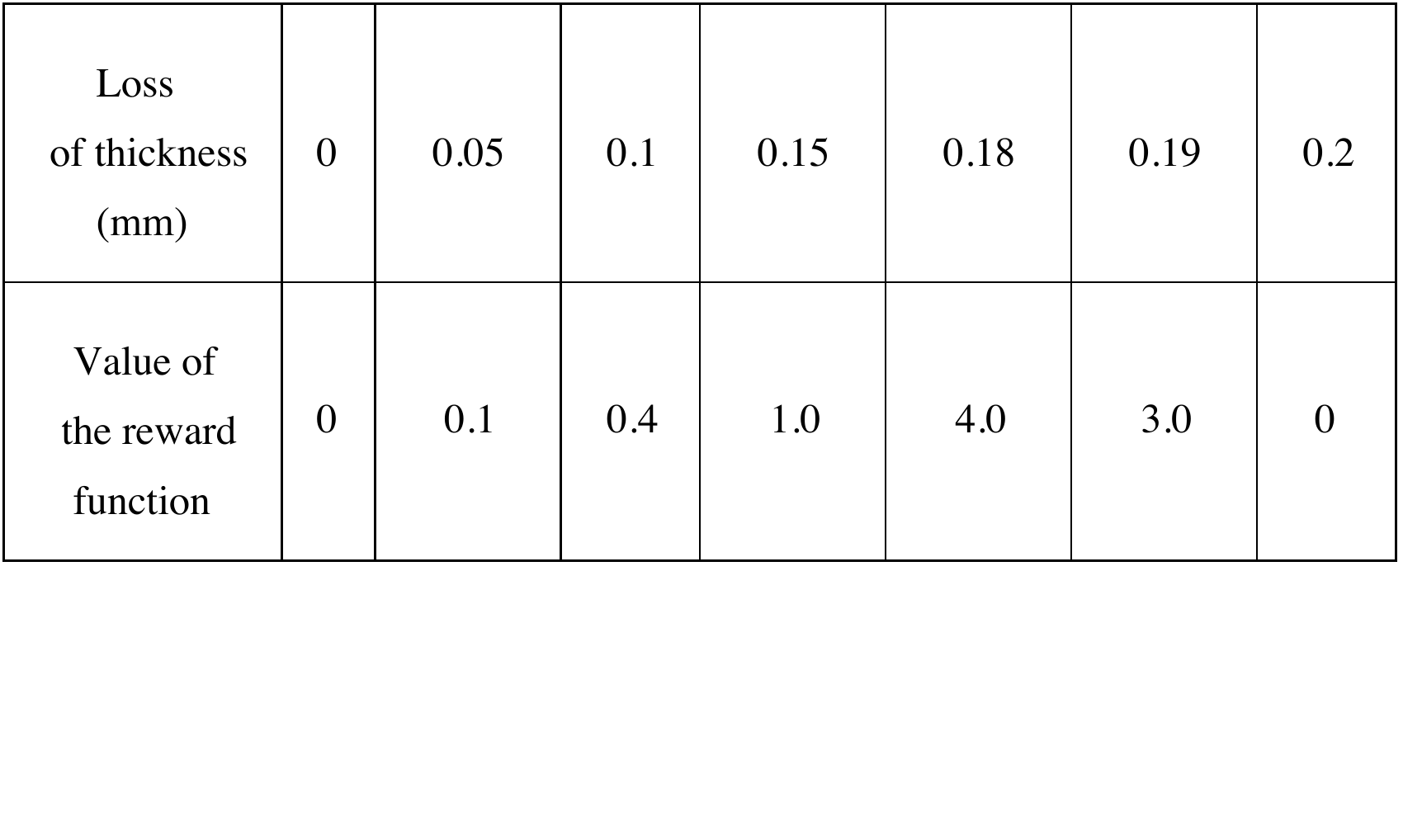}}
\caption{Graphical representation and definition of the cost function as a function of 
the thickness loss}
\label{figure_3}
\end{figure}

For our example of metallic structure, we choose an arbitrary reward function that depends 
only on the loss of thickness, since this is the critical factor to monitor. Note that we could take into account the other components of our process without any additional difficulty.
The reward function is built to reflect the fact that beyond a loss of thickness of 0.2mm, 
the structure is unusable, so it is too late to perform maintenance. 
Conversely, if the thickness loss is small, such a maintenance is 
unnecessarily costly. We use a piecewise affine function $g$ which values are 
given at the points in the table in Figure~\ref{figure_3}.
As for the choice of the computational horizon $N$, numerical simulations show that over 25 changes of environment, all 
trajectories exceed the critical threshold of $0.2$mm. We will therefore set
the  time horizon to be the 25th jump ($N = 25$).

\section{Numerical procedure}
\label{num}

It is natural to propose an iterative algorithm to calculate an approximation  
of the value function based on a discretization of the operator $L$ defined in equation~(\ref{def L}). This poses 
several problems, related to maximizing continuous functions, the presence of 
the indicator and the presence of conditional expectations. We nevertheless 
managed to overcome these three problems, using the specific properties of PDMP's, 
and in particular the fact that the operator $L$ depends only on the Markov chain 
$(Z_n, S_n)$. Our algorithm for calculating the value function is divided into 
three stages described below: a quantization of the Markov 
chain $(Z_n, S_n)$, a path-adapted time discretization between jumps, and 
finally a recursive computation of the value function $v$. Then, 
the calculation of quasi-optimal stopping time only uses comparisons of 
quantities already calculated in the approximation of the value function, 
which makes this technique particularly attractive, see \cite{saporta10a} for more mathematical details.

\subsection{Quantization}
The goal of the quantization step is to replace the continuous state space 
Markov chain $(Z_n, S_n)$ 
by a  discrete state space chain 
$(\hat{Z}_n, \hat{S}_n)$. The quantization algorithm is described in 
details in e.g. \cite{pages98a} \cite{pages05a} \cite{pages04a} or \cite{pages04b}.
The principle is to obtain a finite grid adapted to the distribution of the 
random variable, rather than building an arbitrary regular grid. We discretize random variables rather than the state space, the idea is to put more points in the areas oh high density of the random variable. The quantization algorithm is based 
on Monte Carlo simulations combined with a stochastic gradient method. 
It provides $N+1$ grids
$\Gamma_n, ~ 0\leq n \leq N $ of dimension $d+2$, one for each couple $(Z_n,S_n)$, with $K$ points in each grid. The algorithm also
provide weights for the grid points and probability transition between two points of two
consecutive grids. 

\bigskip

We note $p_n$ the projection to the nearest neighbor (for the Euclidean norm) 
from $\mathbb{R}^{d+2}$ onto $\Gamma_n$. The approximation of the Markov chain 
$(Z_n, S_n)$ is constructed as follows:
\begin{displaymath} 
  (\hat{Z}_n, \hat{S}_n) = p_n (Z_n, S_n).
\end{displaymath} 
Note that $\hat{Z}_n$ and $\hat{S}_n$ depend on both $Z_n$ and $S_n$. The quantization theory ensures that the $L^2$ norm of the 
distance between $(\hat{Z}_n, \hat{S}_n)$ and $(Z_n, S_n)$ tends to 0 as the 
number of points $K$ in the  quantization grids tends to infinity, see  \cite{pages04a}.

\bigskip

\begin{figure}[ht]
\centering
\subfigure[Standard normal density in 2D]
{\includegraphics[width=0.49\linewidth]{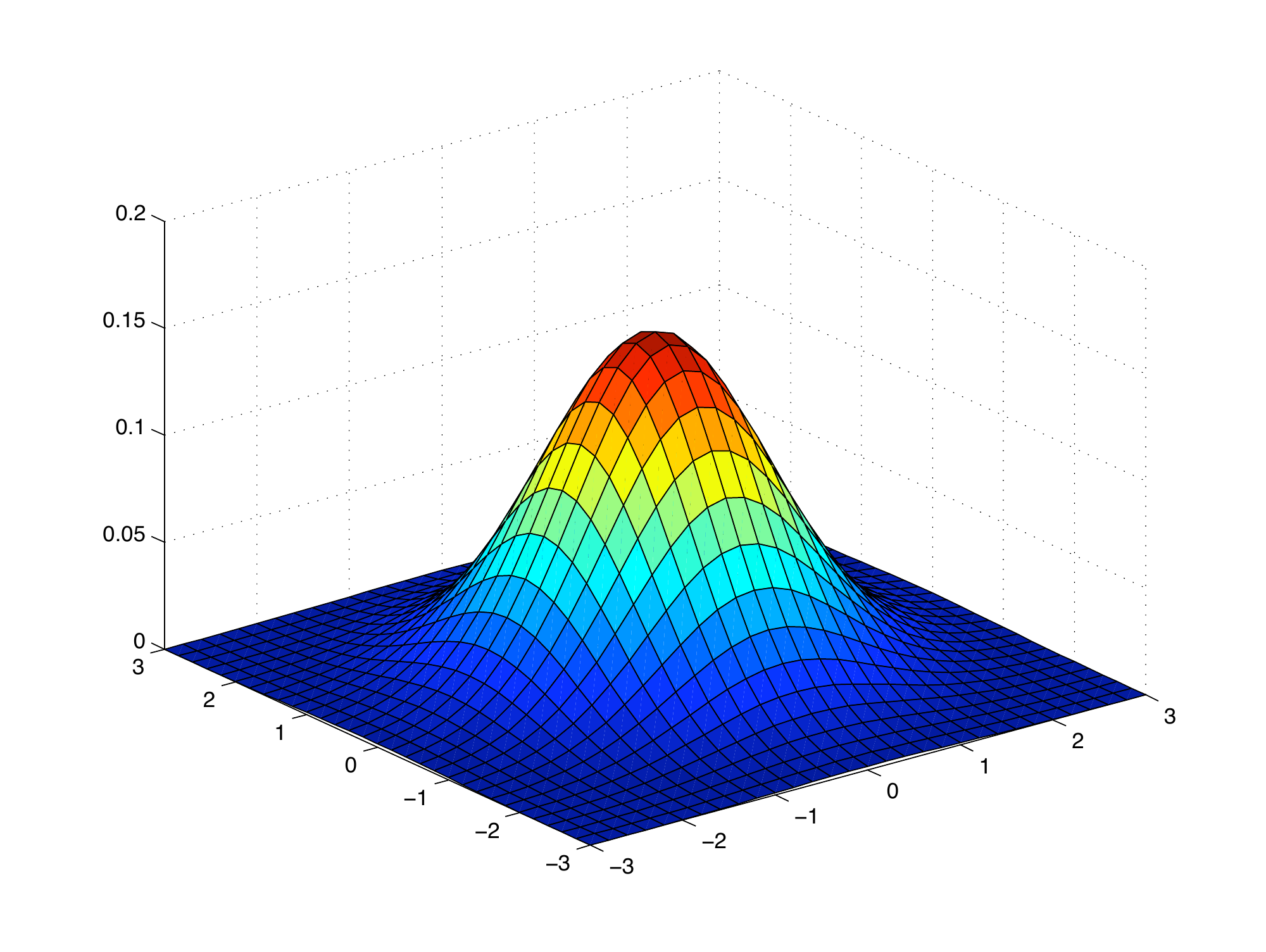}} 
\subfigure[Quantization grid]
{\includegraphics[width=0.49\linewidth]{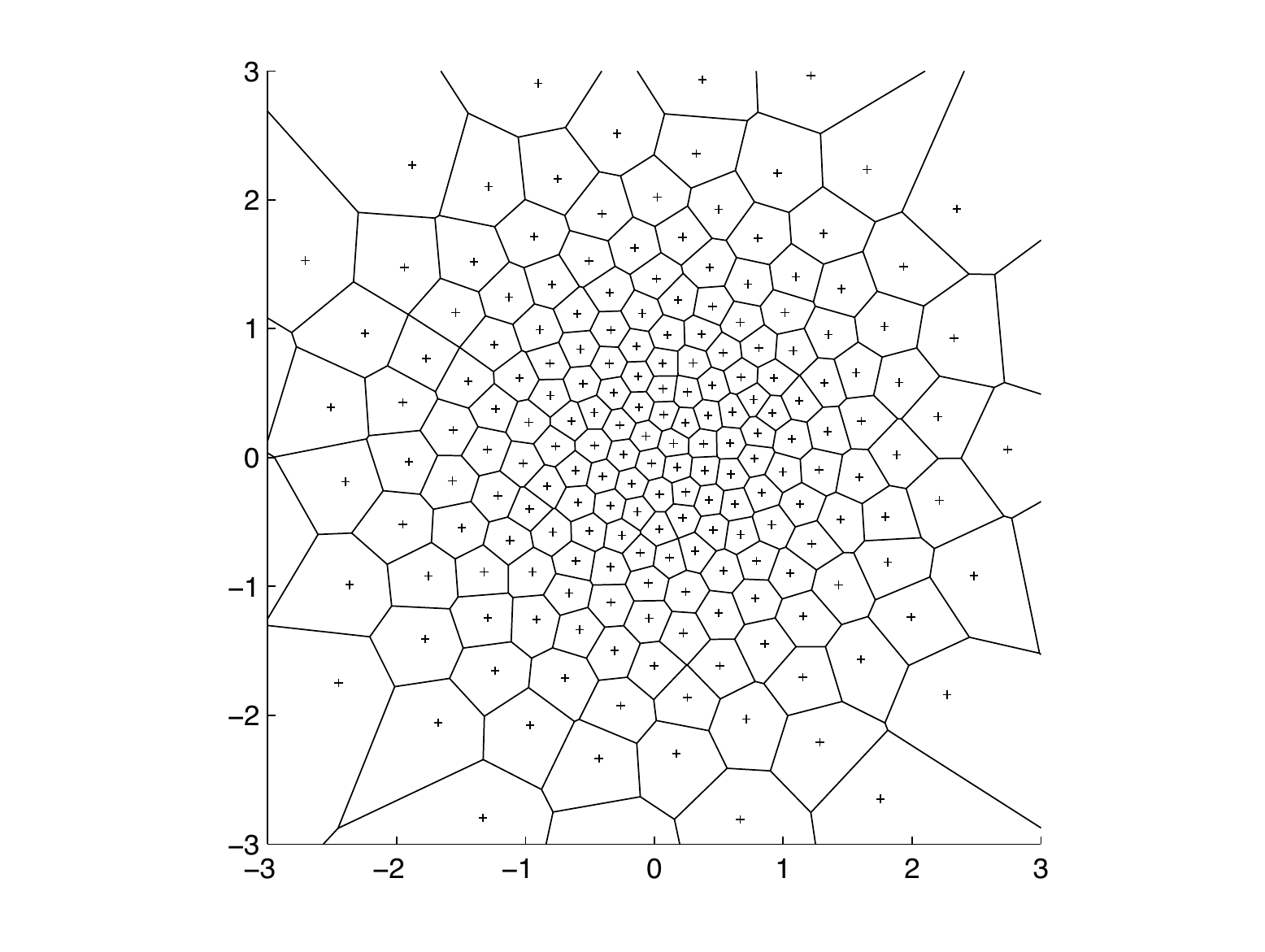}} 
\caption{Example of quantization grid for a normal distribution}
\label{figure_4}
\end{figure}

It should be noted that when the dimension of $Z$ is large, $N$ is large and 
we want to obtain grids with a large number $K$ of points, the quantization algorithm 
can be time-consuming. However, we can make this grids calculation in advance 
and store them. They depend only on the distribution of the process, and not on the cost 
function. 
Figure \ref{figure_4} gives an example of quantization grid for the standard
normal distribution in two dimensions. It illustrates that the quantization algorithm puts more 
points in areas of 
high density.

\subsection{Time discretization}
We now wish to replace the continuous maximization of the operator $L$ by a 
finite maximization, that is to say that we must discretize the time 
intervals $[0, t^\ast(z)]$ for each $z$ in the quantization grids. For this, we choose a time step $\Delta <t^\ast(z)$ 
(which may depend on $z$) and we construct the grids 
$G(z) = \{t_1,\cdots, t_{n (z)}\}$ defined by
\begin{itemize} 
\item $n(z)$ is the integer part minus 1 of $t^\ast(z)/\Delta$,
\item for $1 \leq i \leq n (z)$, $t_i = i\Delta$.
\end{itemize} 
We obtain grids that not only do not contain $t^\ast(z)$, but in addition, 
their maximum is strictly less than $t^\ast(z) - \Delta$, which is a crucial 
property to derive error bounds for our algorithm, see \cite{saporta10a}. Note also that we only need a finite number of grids $G(z)$, corresponding 
to the $z$ in the quantization grids  
$(\Gamma_n)\, 0 \leq n \leq N$. Calculation of this time grids can still be made 
in advance. Another solution is to store only $\Delta$ and  $n(z)$ which are 
sufficient to reconstruct the grids.

\bigskip

In practice, we choose a $\Delta$ that does not depend on $z$. To ensure that 
we have no empty grid, we first calculate the minimum of $t^\ast(z)$ on all grids 
of quantization, then we choose a $\Delta$ adapted to this value.

\subsection{Approximate calculation of the value function}
We now have all the tools to provide an approximation of the operator $L$. For 
each $1\leq n \leq N$, and for all $z$ in the quantization grid at time $n-1$, 
we set
\begin{eqnarray*} 
 \hat{L}_n (w,g)(z) 
&\equiv& \max_{u\leq G(z) }
\Big \{ 
E  \Big [ w(\hat{Z}_{n-1})1_{\hat{S}_{n}<u\wedge t^\ast(z)} + g(\Phi(\hat{Z}_{n-1},u))1_{\hat{S}_n\geq u \wedge t^\ast(z)} 
  |\hat{Z}_{n-1}=z \Big ] \Big \}\\
&& \qquad \vee E \left [ w(\hat{Z}_n)|\hat{Z}_{n-1}=z \right ].
\end{eqnarray*} 
Note that because we have different quantized approximations at each time step, we also have different discretizations of operator $L$ at each time step. We then construct an approximation of the value function by backward iterations of the $ \hat{L}_n$:
\begin{displaymath} 
\left \{
\begin{array} {lcl}
\hat{v}_N & = & g, \\
\hat{v}_{n-1} (\hat{Z}_{n-1})& = & \hat{L}_n(\hat{v}_{n},g)(\hat{Z}_{n-1}), 
\quad 1\leq n\leq N.
\end{array} 
\right .
\end{displaymath} 
Then we take $\hat{v}_0(\hat{Z}_0)=\hat{v}_0(z)$ as an approximation of the 
value function $v$ at the starting point $z$ of the PDMP. It should be noted that the 
conditional expectations taken with respect to a process with discrete state space 
are actually finite weighted sums.

\begin{theorem} 
Under assumptions of Lipschitz regularity of the cost function $g$ and local 
characteristics $(\Phi,\lambda,Q)$ of the PDMP, the approximation error in the 
calculation of the value function is
\begin{displaymath} 
\displaystyle 
||\hat{v}_0(z) - v_0(z)||_2 \leq C \sqrt{E Q}
\end{displaymath} 
where $C$ is an explicit constant which depends on the cost function and 
local characteristics of the PDMP, and $EQ$ is the quantization error.
\end{theorem} 

Since the quantization error tends to 0 when  the number of points in the 
quantization grid increases, this result shows the convergence of our 
procedure. Here, the order of magnitude as the 
square root of the quantization error
is due to the presence of indicator functions, which slow convergence 
because of 
their irregularity. To get around the fact that these functions are not 
continuous, we use the fact that the sets where they are actually discontinuous 
are of very low probability. The precise statement of this theorem and its proof can be found in 
\cite{saporta10a}.

\subsection{Calculation of a quasi-optimal stopping time}
We have also implemented a method to compute an $\epsilon$-optimal stopping 
time. The discretization is much more complicated and subtle 
than that of operator $L$, because we need both to use the true Markov chain 
$(Z_n, S_n)$ and its quantized version $(\hat{Z}_n, \hat{S}_n)$. 
The principle is as follows:
\begin{itemize} 
\item At time $0$, with the values $Z_0 = z$ and $S_0 = 0$, we calculate  a first 
date $R_1$ which depends on $Z_0$, $S_0$ and on the value that has realized 
the maximum in the calculation of $\hat{L}_1(\hat{v}_1,g)$.
\item We then allow the process to run normally until the time 
$R_1 \wedge T_1$, that is the minimum between this computed time $R_1$ and the first change of environment. If $R_1$ comes first, it is the date of near-optimal maintenance, 
if $T_1$ comes first, we reset the calculation.
\item At time $T_1$, with the values of $Z_1$ and $S_1$, we calculate  the second 
date $R_2$ which depends on $Z_1$ and $S_1$ and on the  
the value that has realized the maximum in the calculation of 
$\hat{L}_2(\hat{v}_2,g)$.
\item We then allow the process to run normally until the time 
$(T_1 + R_2) \wedge T_2$, that is the minimum between the computed remaining time $R_2$ and the next change of environment. If $T_1 + R_2$ comes first, it is the date of near-optimal
 maintenance, if $T_2$ comes first, we reset the calculation, 
and so on until the $N$th  jump time where maintenance will be performed 
if it has not occurred before.
\end{itemize} 
We have also proved the quality of this approximation by comparing the 
expectation of the cost function of the process stopped by the above strategy 
to the true value function. This result, its proof and the precise construction 
of our stopping time procedure can be found in \cite{saporta10a}.

\bigskip

This stopping strategy is interesting for several reasons. First, this is a 
real 
stopping time for the original PDMP which is a very strong result. Second, it requires no additional 
computation compared to those made to approximate the value function. 
This procedure can be easily performed in real time, and only requires an observation of the process at the times of change of environment, which is exactly the available inspection data for our metallic structure. Moreover, even if the 
original problem is an optimization {\em on average}, this stopping rule is 
path-wise and is updated when new data arrive on the history of the process at each change of environment.
Finally, as our 
stopping procedure is of the form {\em intervene at such date if no change of environment
has occurred in the meantime},  it allows in some measure to have  
maintenance scheduled 
in advance, 
In particular, our procedure ensures that there will be no need to perform 
maintenance before a given date, which is crucial for our example as a submarine in operation should not be stopped  at short notice.

\section{Numerical results}
\label{results}

We have implemented this procedure for the optimization of the maintenance of 
the metallic structure described in section~\ref{modeling}. With our choice of reward function, it is easy to see that the true value function at $z=0$ is 4, which is 
the maximum of the reward function $g$, and an optimal stopping time is the first moment 
when the loss reaches 0.18 mm thick (value where $g$ reaches its maximum). This 
is because the cost function only depends on the thickness loss, which  
evolves continuously increasingly over time. However, our numerical procedure is valid for any sufficiently regular reward function, and we shall not use the knowledge of the true value function or optimal stopping time in our numerical procedure. Besides, we recall that the thickness loss is not measured continuously.

\bigskip

While running the algorithm described in the previous section, e encountered an unexpected difficulty for the construction of the 
quantization 
grids. 
Indeed, the scales of the different variables of the problem are
radically different: from about $10^{-6}$  for $\rho$ to $10^5$ for the average time 
spent in environment 2. This poses a problem in the classical quantization 
algorithm as searching the nearest neighbor and gradient calculations 
are done in Euclidean norm, regardless of the magnitudes of the components. 
\begin{figure}[ht]
\centering
\subfigure[Classical algorithm]
{\includegraphics[width=0.49\linewidth]{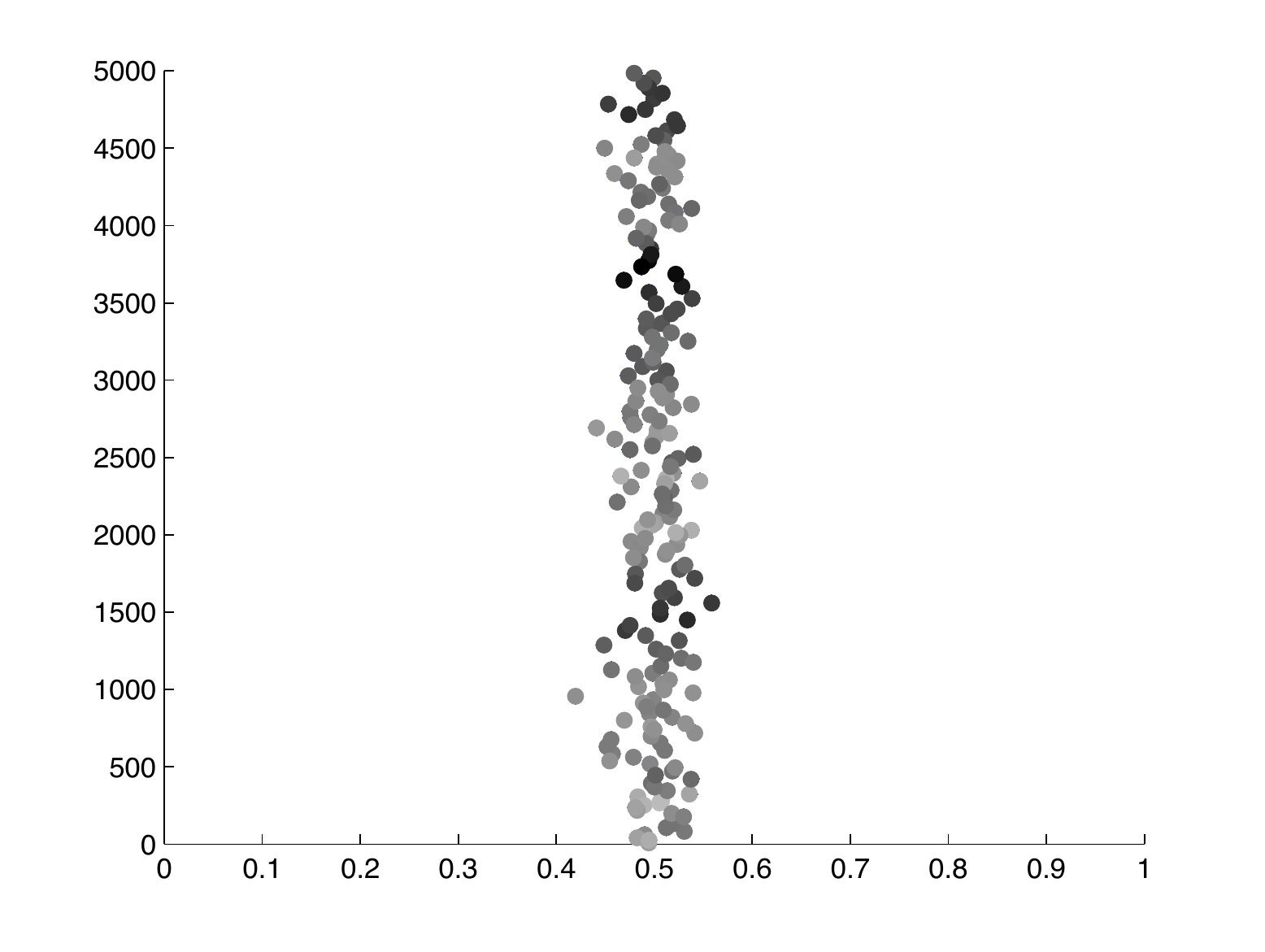}} 
\subfigure[Algorithm with weighted Euclidean norm]
{\includegraphics[width=0.49\linewidth]{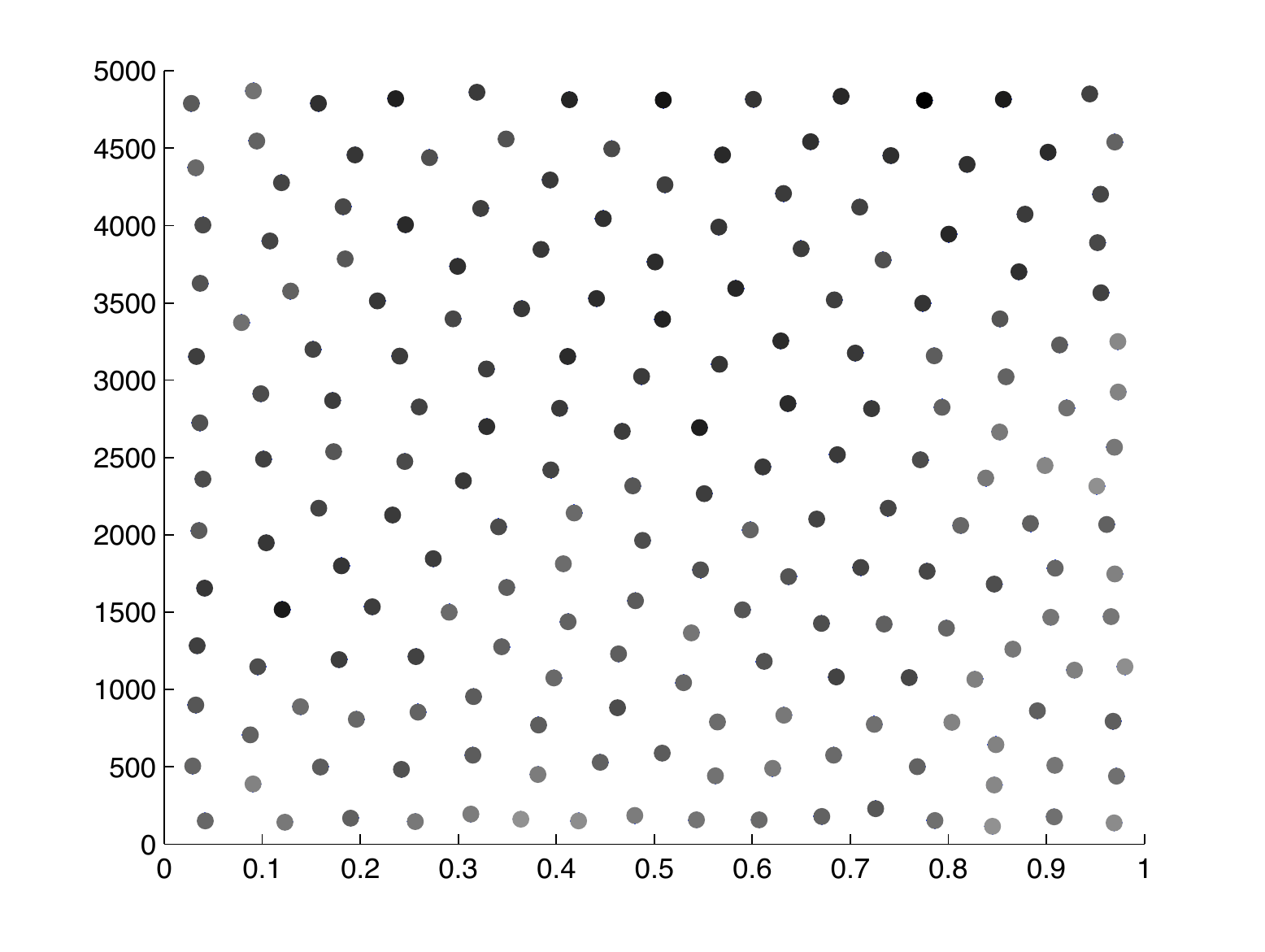}} 
\caption{Quantization grids for a uniform distribution on $[0,1]\times [0, 5000]$} 
\label{figure_5}
\end{figure}
Figure \ref{figure_5} illustrates this problem by presenting two examples of 
quantization grids for a uniform distribution on $[0,1]\times [0,5000]$. The left image 
shows the result obtained by the conventional algorithm, the right one is obtained 
by weighting the Euclidean norm to renormalize each variable on the same scale. 
It is clear from this example that the conventional method is not satisfactory, 
because the grid obtained is far from uniform. This defect is corrected by a 
renormalization of the variables. We therefore used a weighted Euclidean 
norm to 
quantify the Markov chain associated with our degradation process. 

\bigskip

\begin{figure}[ht]
\centering
\subfigure[ Environment 2, time $T_1$]
{\includegraphics[width=0.32\linewidth]{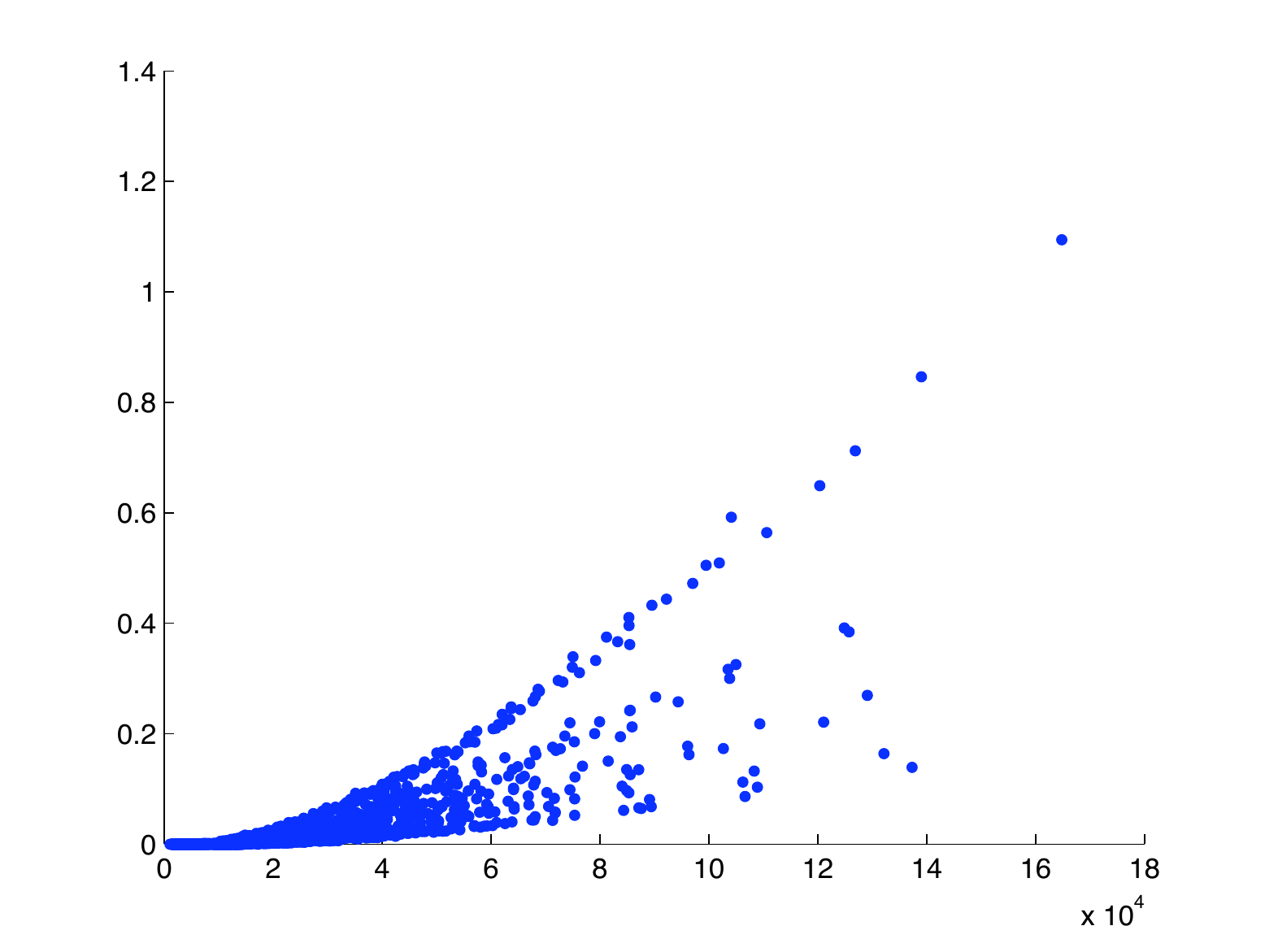}}
\subfigure[ Environment 3, time $T_2$]
{\includegraphics[width=0.32\linewidth]{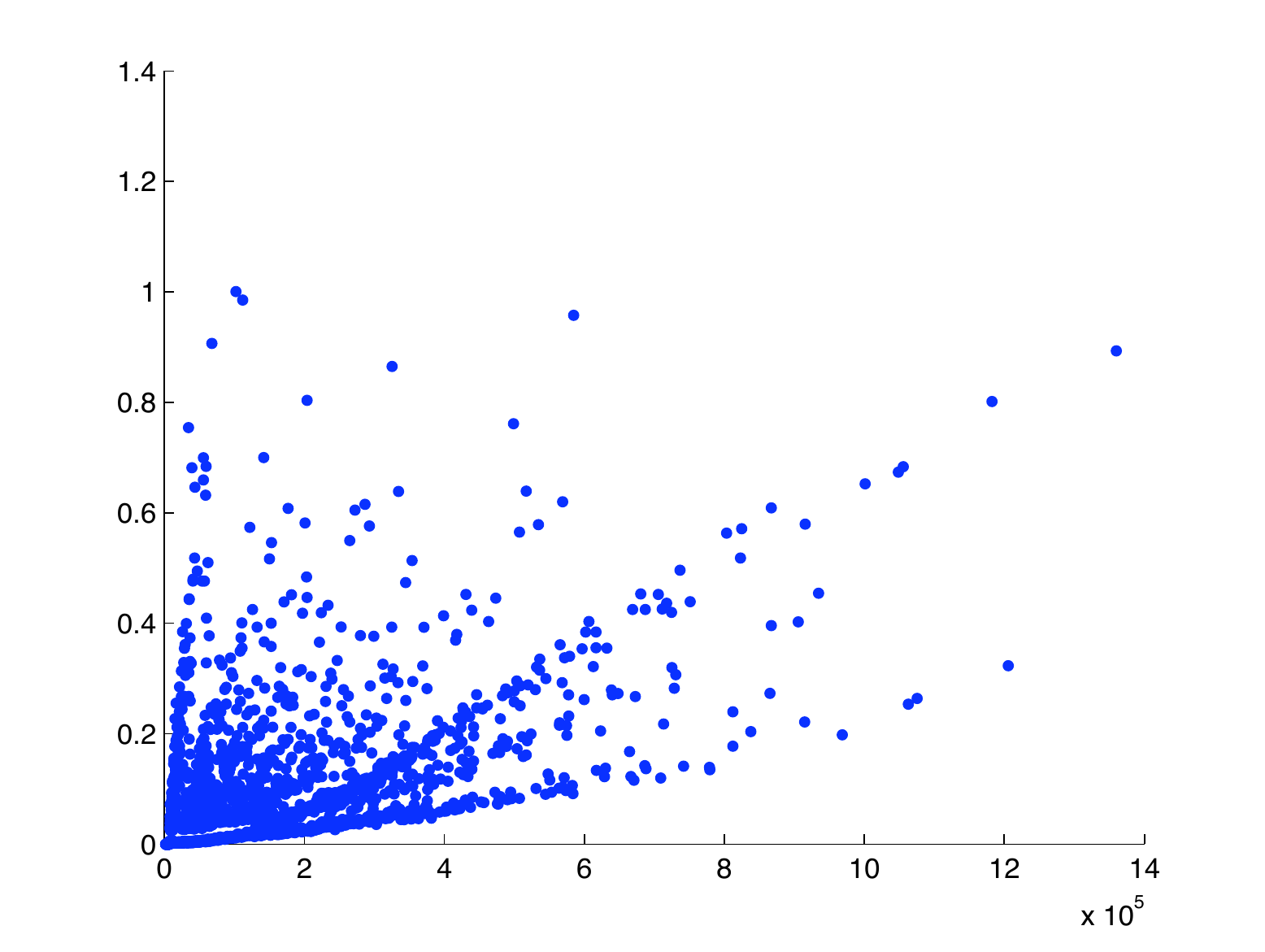}}
\subfigure[ Environment 2, time $T_{10}$]
{\includegraphics[width=0.32\linewidth]{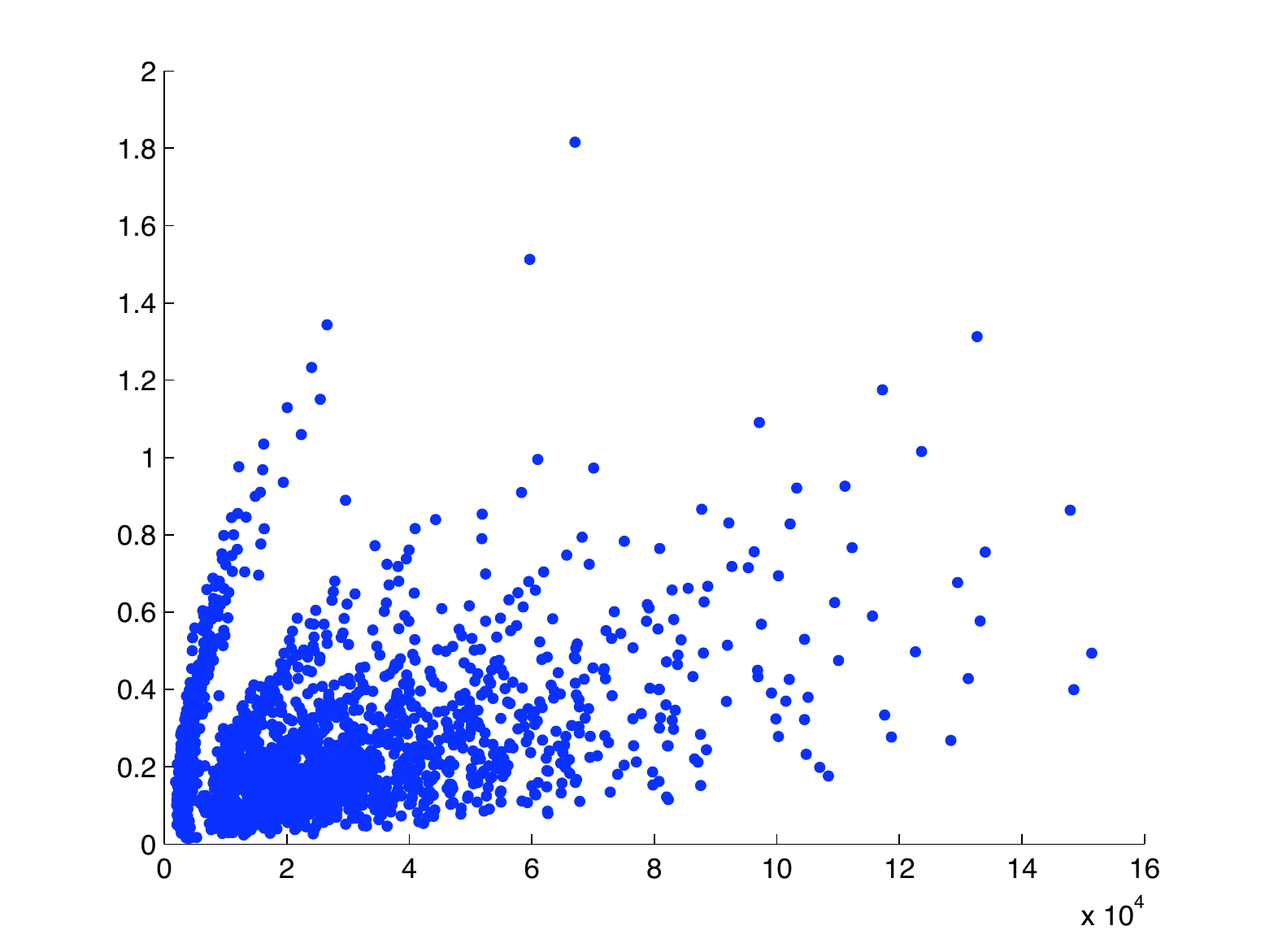}}
\subfigure[ Environment 1, time $T_{15}$]
{\includegraphics[width=0.32\linewidth]{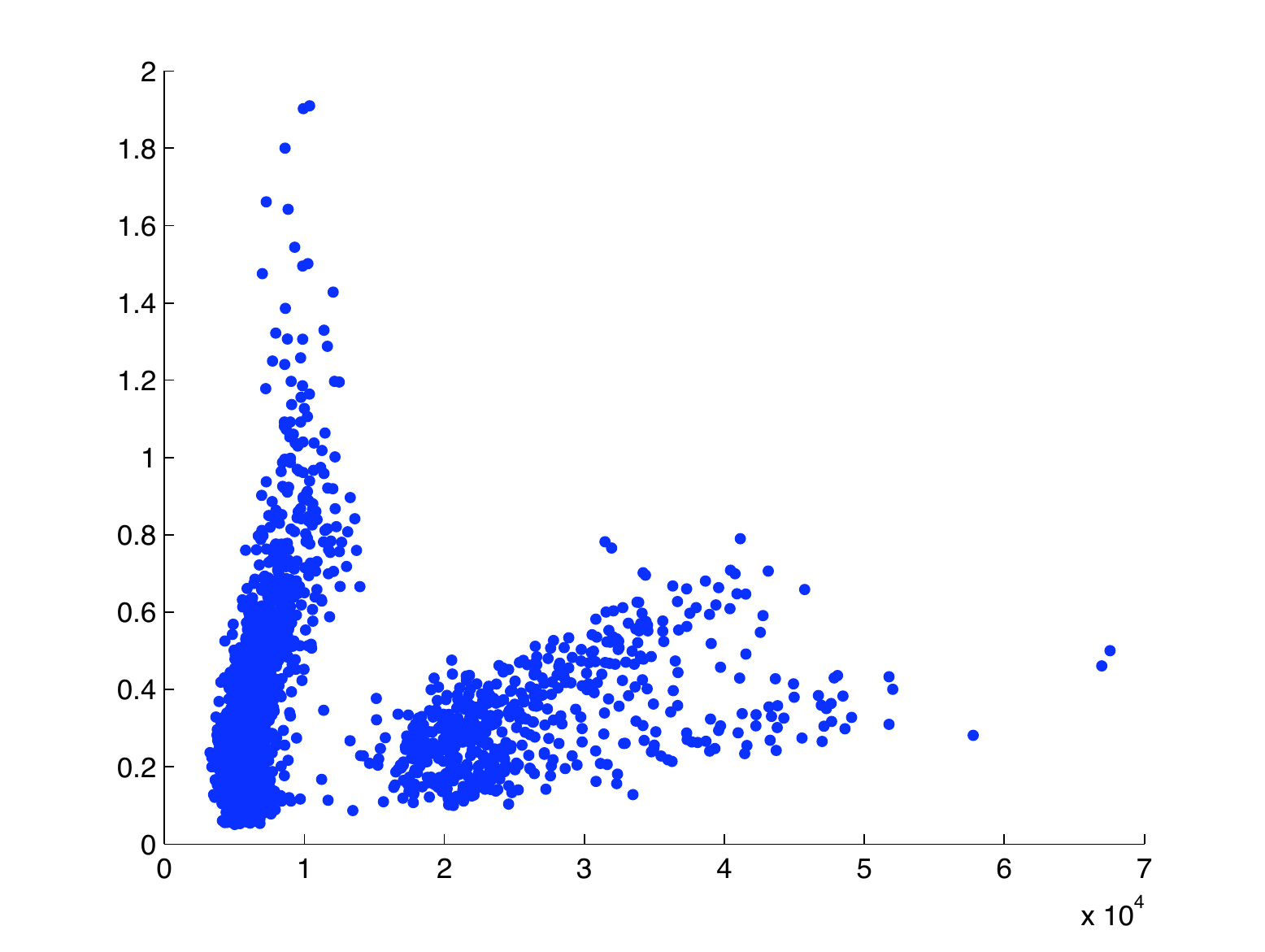}}
\subfigure[ Environment 2, time $T_{19}$]
{\includegraphics[width=0.32\linewidth]{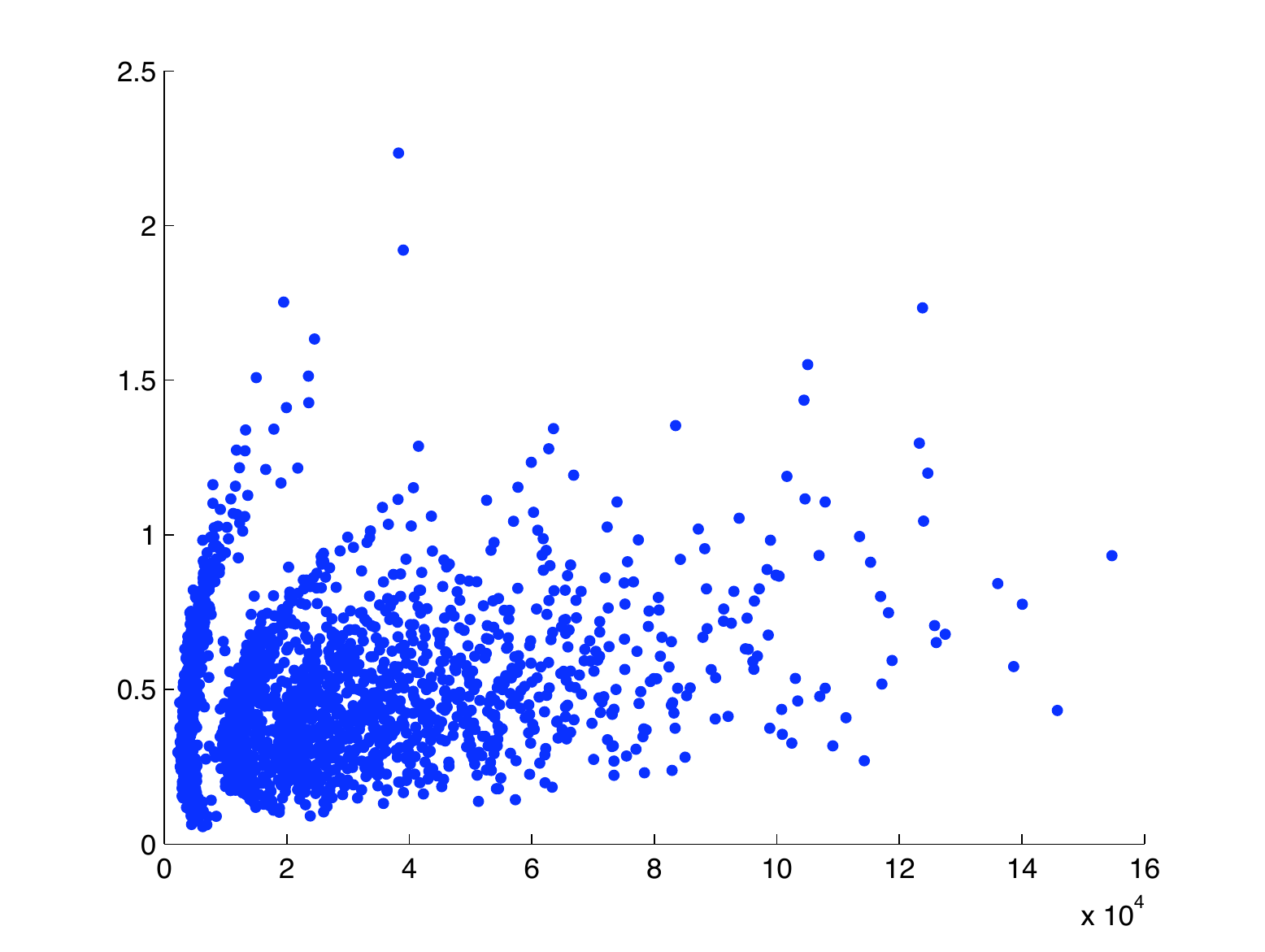}}
\subfigure[ Environment 2,time $T_{25}$]
{\includegraphics[width=0.32\linewidth]{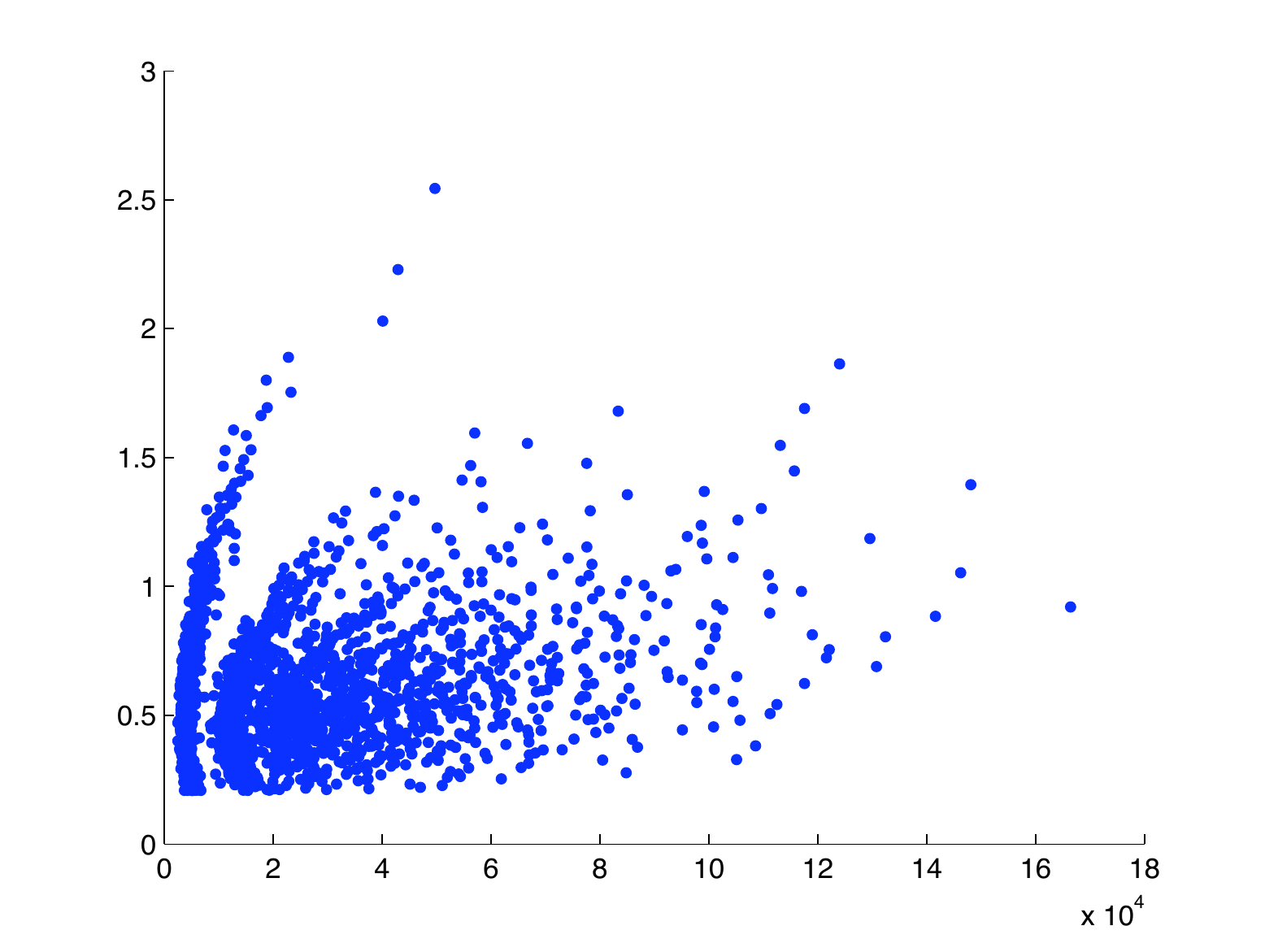}}
\caption{Quantization grids with 2000 points for the inter-jump time (abscissa) 
and the thickness loss (ordinate). The scale changes for each graph.} 
\label{figure_6}
\end{figure}

Figure \ref{figure_6} 
shows some projections of the quantization 
grids with 2000 points 
that we obtained. The times are chosen in order to to illustrate the random and 
irregular nature of the grids, they are custom built to best approach the distribution 
of the degradation process.

\bigskip

\begin{figure}[ht]
\centering
\subfigure
{\includegraphics[width=0.49\linewidth]{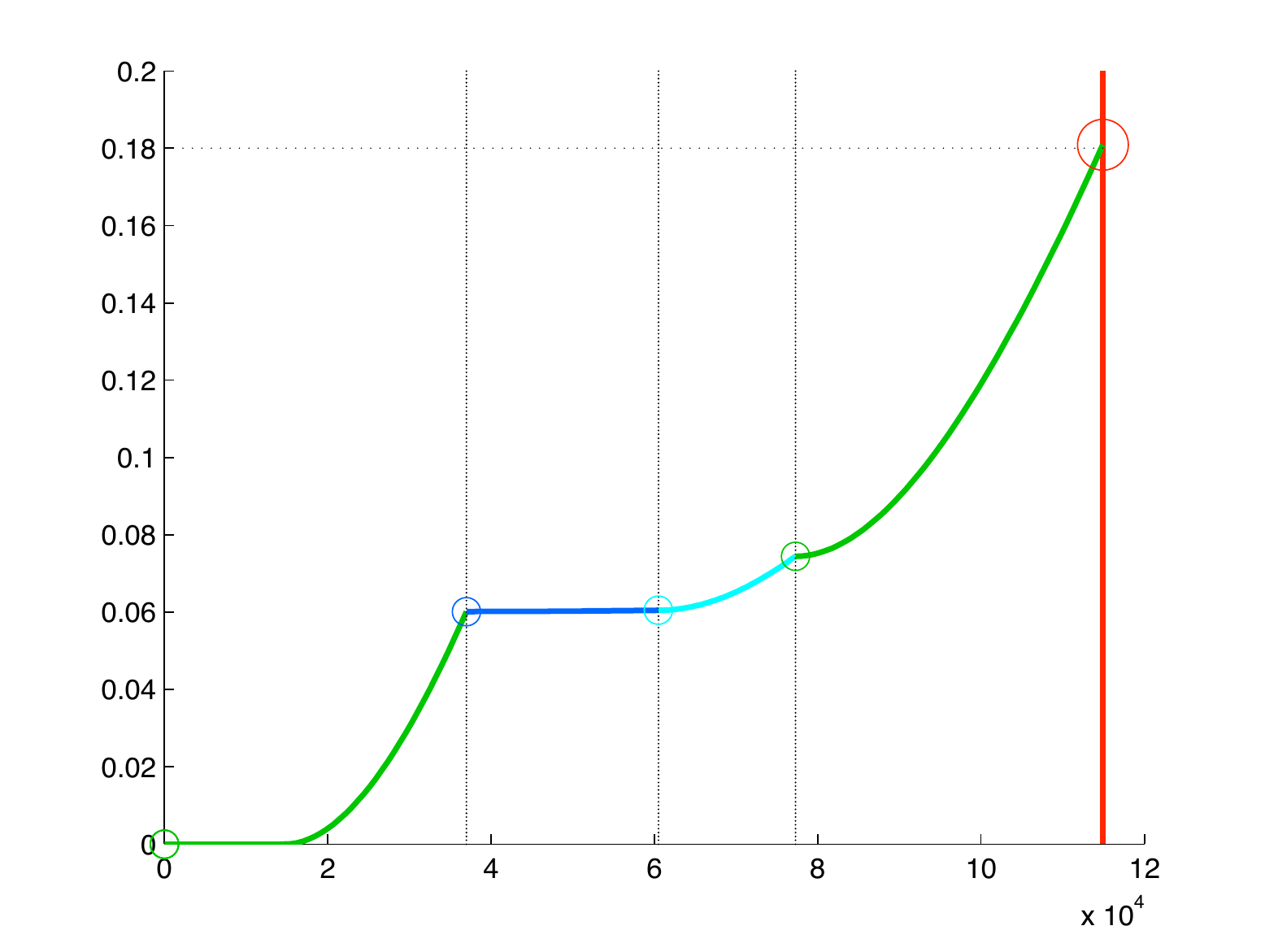}}
\subfigure
{\includegraphics[width=0.49\linewidth]{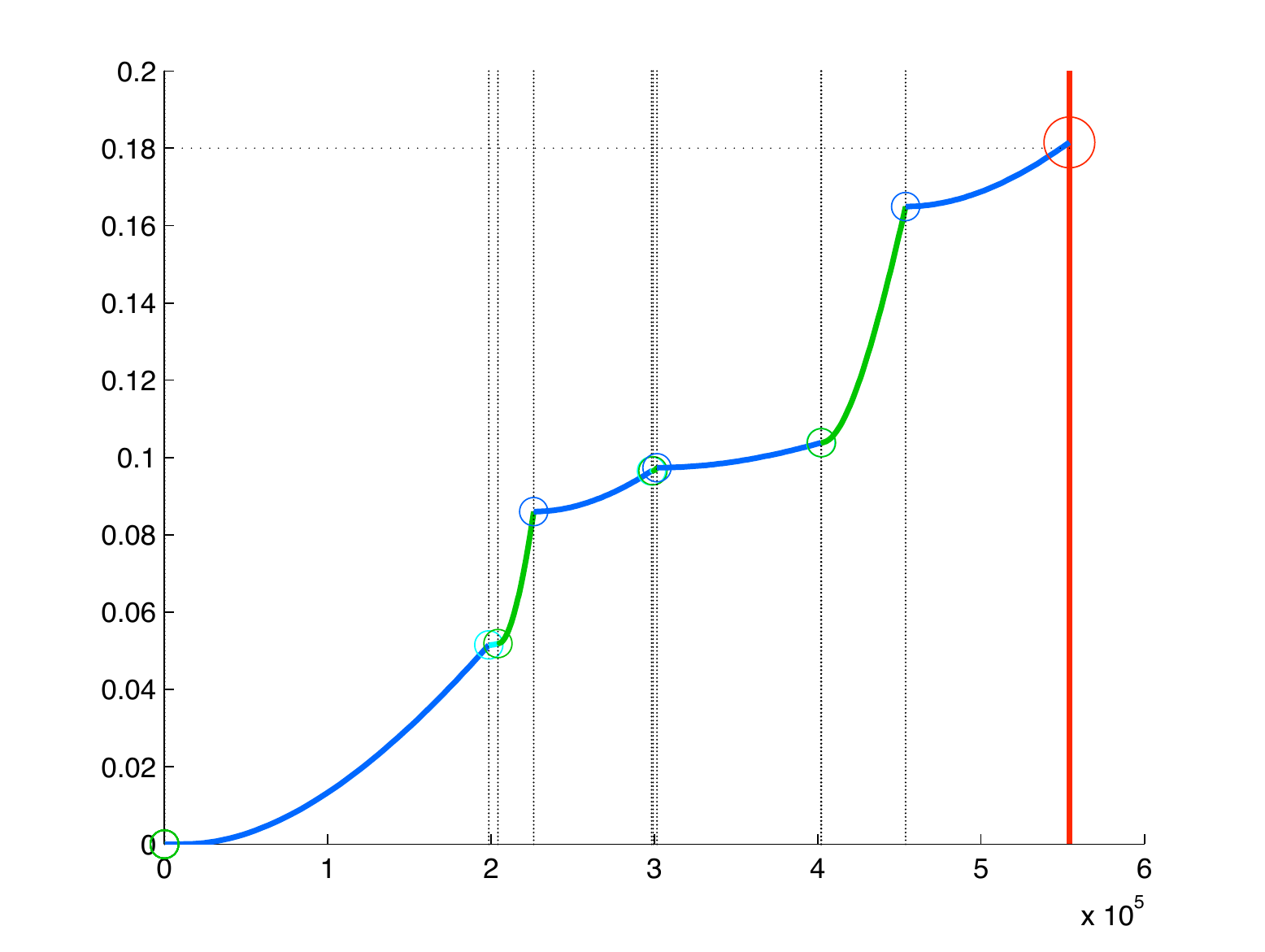}}
\caption{Examples of stopped trajectories with the optimal maintenance time 
calculated by the algorithm.} 
\label{figure_7}
\end{figure}

Figure \ref{figure_7} 
shows two examples of computation of the quasi optimal maintenance time on two specific simulated trajectories. The thick vertical line 
represents the moment provided by the algorithm to perform maintenance. 
The other vertical lines materialize the moments of change of  environment, 
the horizontal 
dotted line the theoretical optimum. In both examples, we stop at a value very 
close to the optimum value. In addition, the intervention did take place before the critical 
threshold of 0.2mm.

\bigskip

We calculated an approximate value function $v$ in two ways. The first one is the 
direct method obtained by the algorithm described above. The second one is obtained 
by Monte Carlo simulation using the quasi-optimal stopping time provided by 
our procedure. The numerical results we obtained are summarized in 
Table \ref{table_1}.
\begin{table}[ht]
\begin{center} 
\begin{tabular} {|c|c|c|} \hline
Number of points             & Approximation of the            & Approximation of the value         \\
in the quantization        & value function by the           &  function by Monte Carlo  with the \\
grids                        & direct algorithm                       &  quasi-optimal stopping time       \\ \hline
10   &   2.48   & 0.94 \\ \hline
50   &   2.70   & 1.84 \\ \hline
100  &  2.94	& 2.10 \\ \hline
200  & 3.09	& 2.63 \\ \hline
500  & 3.39     & 3.15 \\ \hline
1000 & 3.56	& 3.43 \\ \hline
2000 & 3.70	& 3.60 \\ \hline
5000 & 3.82	& 3.73 \\ \hline
8000 & 3.86     & 3.75 \\ \hline
\end{tabular} 
\caption{Numerical results for the calculation of the value function.}
\label{table_1}
\end{center} 
\end{table} 
We see as expected, that the greater the number of points 
in the quantization grid, the better our approximation becomes. 
Furthermore, the specific form of this cost function $g$ indicates that 
at the threshold of 1, the intervention takes place between 0.15 and 0.2mm, 
and when the threshold increases, this range is narrowed. We can therefore 
state that our approximation is good even for low numbers of grid points. 
The last column of the table also shows the validity of our stopping rule. 
It should be noted here that this rule does not use the optimal stopping 
time {\em stop at the first moment when the thickness loss reaches 0.18mm}. 
The method we use is general and implementable even when the 
optimal stopping time is unknown or does not exist.

\bigskip

\begin{figure}[ht]
\centering
\subfigure[Histogram of 100000 values of the optimal maintenance time expressed 
in years.]
{\includegraphics[width=0.49\linewidth]{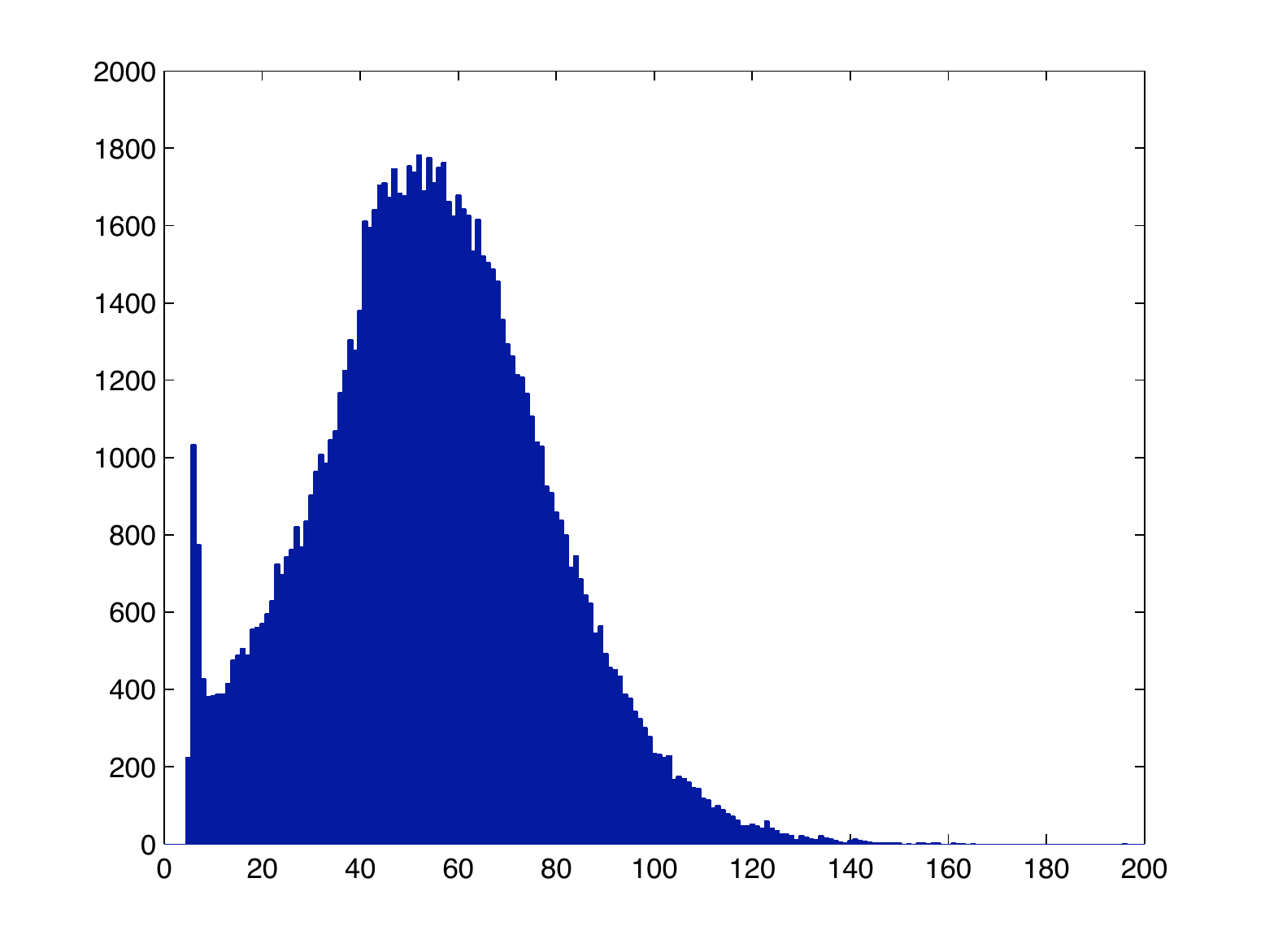}} 
\subfigure[Quantiles.]
{\includegraphics[width=0.49\linewidth]{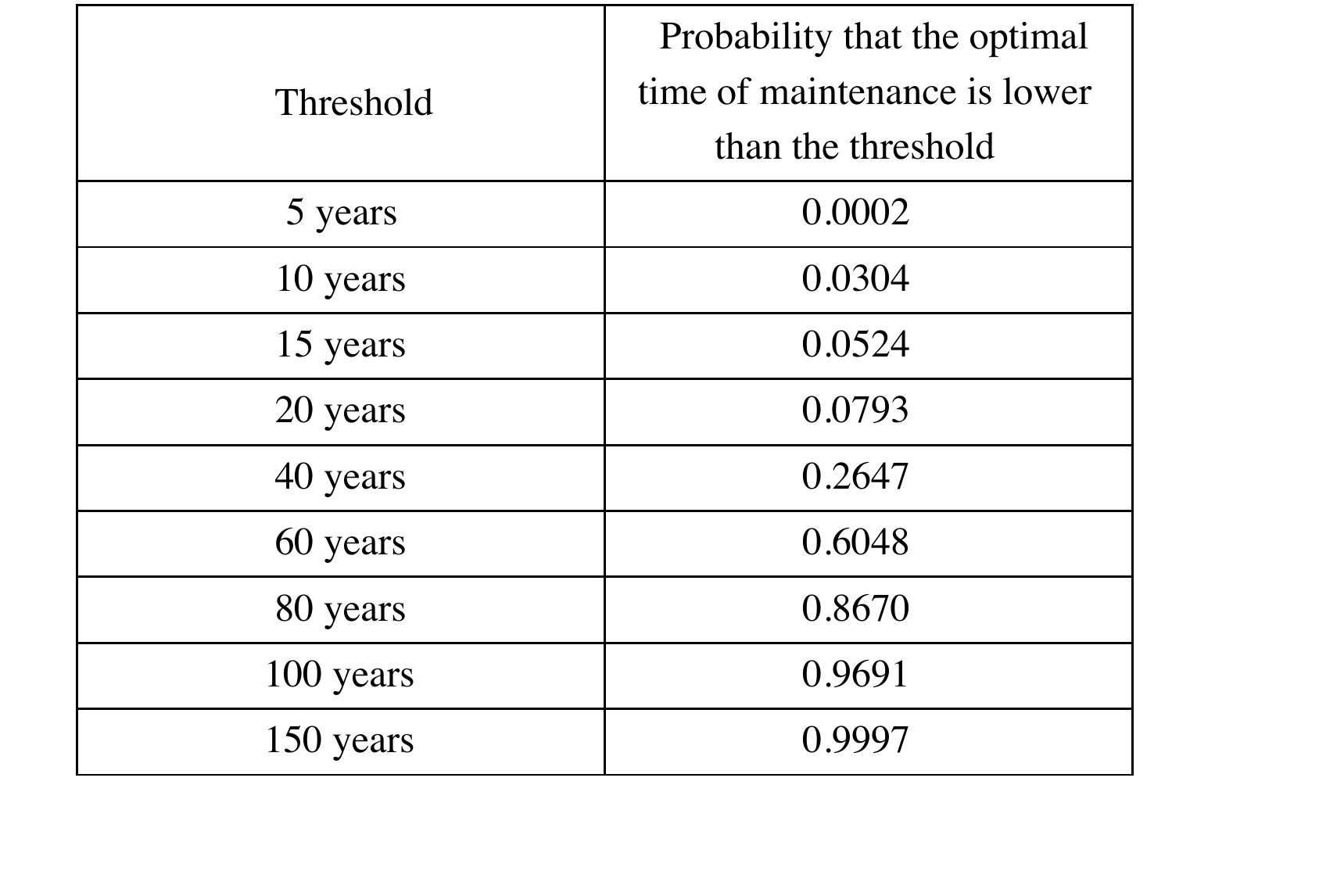}} 
\caption{Distribution and quantiles of the quasi-optimal stopping time.} 
\label{figure_8}
\end{figure}

Moreover, we can also construct a histogram (Figure \ref{figure_8}) 
of the values of our stopping 
time, that is to say, a histogram of the values of effective moments of 
maintenance. We can also estimate the probability that this moment is 
below certain thresholds. These results are interesting for Astrium in the design phase of the 
structure to optimize margins from the specifications and to consolidate 
the design margins available. Thus, we can justify that with a given probability 
no maintenance will be required before the termination date of the contract.

\section{Conclusion}
\label{ccl}

We have applied the numerical method described in \cite{saporta10a} on a practical industrial example to approximate 
the value function of the optimal stopping problem and a quasi-optimal 
stopping time for a piecewise-deterministic Markov process, that is the quasi optimal maintenance date for our structure.
The quantization method we propose can sometimes be costly in computing time, 
but has a very interesting property: it can be calculated off-line. Moreover 
it depends only on the evolutionary characteristics of the model, and not on the 
cost function chosen, or the actual trajectory of the specific process we want to 
monitor.
The calculation of the optimal maintenance time is done in real time.
This method is especially attractive as its application requires knowledge of 
the system state only at moments of change of  environment and not in continous 
time.
The optimal maintenance time is updated at the moments when the system 
switches to another environment and has the form 
{\em intervene at such date if no change of mode takes place in the meantime}, 
which allows to  schedule maintenance services in advance.

\bigskip

We have implemented this method on an example of optimization of the 
maintenance 
of a metallic structure subject to corrosion, and we obtained very satisfactory results,
 very close to theoretical values, despite the relatively large size of the
 problem.
These results are interesting for Astrium in the design phase of the 
structure 
to maximize margins from the specifications and to consolidate the 
avaible dimensional margins. Thus, we propose tools to justify that with a given 
probability that no 
maintenance will be required before the end of the contract.

\bigskip

The application that we have presented here is an example of maintenance 
{\em as good as new} of the system. The next step will be to allow only partial repair 
of the system. The problem will then be to find simultaneously the optimal 
times  
of maintenance and optimal repair levels. Mathematically, it is an impulse 
control 
problem, which complexity exceeds widely that of the optimal stopping. Here 
again, 
the problem is solved theoretically for PDMP, but there is no practical 
numerical 
method for these processes in the literature. We now work in this direction 
and we hope to be able to extend the results presented above.

\subsubsection*{\underline{Acknowledgement}}
This work was partially funded by the ARPEGE program of National Agency for Research (ANR), 
project FauToCoES, ANR-09-004-SEGI.
\bibliographystyle{plain}
{\small
\bibliography{../../bibtex/zhang}
}
\end{document}